\documentclass[sn-vancouver]{sn-jnl}

\jyear{2021}%

\usepackage{bm}
\newcommand{\N}{\mathcal{N}}
\newcommand{\bt}[1]{{\bf #1}}               

\graphicspath{{figures/}}

\theoremstyle{thmstyleone}%
\theoremstyle{thmstyletwo}%

\theoremstyle{thmstylethree}%

\begin{document}

\title[Quantum diamond magnetometry for navigation in GNSS denied environments]{Quantum diamond magnetometry for navigation in GNSS denied environments}


\author*[1]{Xuezhi Wang}\email{xuezhi.wang@rmit.edu.au}
\author[1]{Wenchao Li}\email{wenchao.li@rmit.edu.au}
\author[1,2]{Bill Moran}\email{wmoran@unimelb.edu.au}
\author[3]{Brant C. Gibson}\email{brant.gibson@rmit.edu.au}
\author[4] {Liam T. Hall}\email{liam.hall@unimelb.edu.au}
\author[5] {David Simpson}\email{simd@unimelb.edu.au}
\author[1,2]{Allison N. Kealy}\email{allisonkealy@gmail.com}
\author[1,3]{Andrew D. Greentree}\email{andrew.greentree@rmit.edu.au}

%
%

\affil*[1]{\orgdiv{School of Science}, \orgname{RMIT University}, \orgaddress{\state{Vic 3000}, \country{Australia}}}

\affil[2]{\orgdiv{Department of Electrical \& Electronic Engineering}, \orgname{University of Melbourne}, \orgaddress{\state{Vic 3053}, \country{Australia}}}

\affil[3]{\orgdiv{ARC Centre of Excellence for Nanoscale BioPhotonics}, \orgname{RMIT University}, \orgaddress{ \state{Victoria 3000}, \country {Australia}}}

\affil[4]{\orgdiv{School of Chemistry}, \orgname{University of Melbourne}, \orgaddress{\state{Vic 3010}, \country{Australia}}}

\affil[5]{\orgdiv{School of Physics}, \orgname{University of Melbourne}, \orgaddress{\state{Vic 3010}, \country{Australia}}}

\abstract{Satellite-based navigation is a transformational technology that underpins almost all aspects of modern life. However, there are environments where global navigation satellite systems (GNSS) are not available, for example undersea or underground, and navigation that is robust to GNSS outages is also required for resilient systems.  In this context, we are exploring the potential for quantum diamond magnetometers to be used as aids for inertial measurement units for navigation in GNSS-denied environments.  We perform simulations of the magnetic field measurements combined with probabilistic data association for data mapping; and probabilistic multiple hypotheses (or Viterbi) map matching filters.  These methods are used to explore the expected navigation errors likely to be present in scalar and vector magnetometry with diamond at the sensitivities available using current and expected near-term devices.}



\maketitle

\section{Introduction}
\label{sec1}

Next-generation quantum technologies are poised to cause disruption to a wide range of sectors \cite{2017:Degen}.  Such potential derives from features including extreme sensitivity, with (in some cases) nano-scale sensor elements, new sensing modalities including entangled sensing \citep{Giovannetti2004}, and the promise of low-Size, Weight and Power (SWAP).  Superconductors represent a mature quantum sensing platform \citep{Porrati2019}, albeit one restricted to cryogenic temperatures, whereas atomic vapour cells \citep{Knappea}  and diamond (the subject of this work) \citep{Taylor2008} are emerging as some of the most promising room-temperature compatible quantum sensors.

Our focus here is on the potential application of diamond containing the negatively-charged nitrogen-vacancy center (NV) for navigation in GNSS denied environments.  NV diamond is promising for such applications due to its ambient temperature and magnetic field operation \citep{Barry2020}, and operation in both projection and vector magnetometry \citep{Schloss2018} modes.  Because NV diamond is a still developing technology, we present out results as a function of sensitivity, noting that  NV magnetometers have demonstrated sensitivity down to 900~fT/$\sqrt{\text{Hz}}$ \citep{Wolf2015}.

In GNSS-denied environments, platform navigation performance is dominated by the accuracy of onboard inertial sensors. Even with high end inertial sensors, which exhibit extremely low bias and drift, it is not possible to avoid the build-up of navigation errors over long time frames~\citep{Titterton2004a}. Removing these accumulated navigation errors is therefore crucial for  navigation accuracy~\citep{Groves2013b}. This removal, or correction, is achieved using one or more aiding sources that provide positional information, i.e. a position fix.

Geophysical map matching is an effective method for localisation and navigation where GNSS is not available; such as underwater, urban, or hostile environments. Such maps may derive from a variety of geophysical data, e.g. gravity, airborne electromagnetics, passive seismic, active seismic, magnetotellurics, magnetic, radiometric, and elevation \citep{Tuohy1996}. Map matching estimates a platform's location by directly matching the geophysical measurement to a point in the map, i.e. the matching occurs explicitly in the map space. The resulting location estimate is then integrated into the navigation system in a similar way to a loosely coupled GNSS/INS system \citep{Titterton2004a}.

Although conceptually simple, map matching with geophysical maps suffers from map measurement ambiguity issues. First, the geophysical measurements themselves are degraded by sensor noise so the measurements will not match the map exactly. Second, the measurements may match multiple points within the map (i.e. non-unique solution) as the map lookup process is a scalar to vector mapping. Third, the locations where the measurements were acquired is of course uncertain. Finally, the map itself suffers from finite spatial and signal resolution.

Various approaches for resolving map measurement ambiguity can be found in the literature. The basic idea is to improve knowledge of the platform location by processing multiple measurements with vehicle trajectory constraints.  For example,  \citep{Tuohy1996} combines multiple geophysical maps (altimetry, magnetic and gravity) with  a multiple hypotheses tracker. In \citep{Kamgar-Parsi1999} a batch of possible platform locations on a map are identified from field measurements by finding the closest contour points using a downhill simplex method to solve the position minimization problem. Following the same idea, an iterated closest contour point (ICCP) algorithm is proposed in \citep{Han2018}, later on improved in~\citep{Wang2022a}. In~\citep{Hostetler1983} a similar technique was used for terrain-aided navigation under the name of multilateral arc matching. Nevertheless, these approaches do not properly model the uncertainties of measurement origin in the statistical sense.

More recently, a probabilistic multiple hypotheses tracking map matching (PMHT-MM) filter was proposed in \citep{Wang2022} for gravitational map matching localisation to aid an inertial navigation system (INS). Unlike the ICCP method, the filter takes both measurement ambiguity and platform kinematic constraints into account in an iterative procedure of forward and backward filtering. The  gravitational map matching localisation is also approached using the Viterbi algorithm under the maximum likelihood criterion in \citep{Li2022}. Both approaches show more accurate, robust and effective performance over the standard ICCP method for  the case of INS aiding using gravitational maps.

Here we consider magnetometery aided inertial navigation with total magnetic intensity (TMI) maps. A probabilistic data association (PDA) method is adopted to address the measurement ambiguity problem, which is a nontrivial problem in  map matching localisation using geophysical data maps. We show that the PDA method provides an effective way to map a field measurement into geolocation, but also enables a quantitative analysis of localisation error with respect to the magnetometer noise levels for a given TMI reference map. Furthermore, we implement a magnetometery aided INS using this  method to determine the relationship between magnetometer noise levels and navigation performance.

This manuscript is organised as follows: We first briefly introduce the mechanism for magnetic field sensing with diamond.  Second, the problem of aided inertial navigation is formulated using Bayesian filtering, and the issues associated with magnetometery map matching are presented in Section~\ref{sec2}. We then describe the probabilistic multiple hypothesis map matching method in Section~\ref{sec3}.  In Section~\ref{sec4}, the performance of the proposed algorithm for aiding of INS using online maps are demonstrated in a realistic navigation scenario without GNSS. 

\section{Magnetometry using nitrogen-vacancy centers in diamond}

Quantum sensing with diamond is typically performed using the negatively charged nitrogen-vacancy (NV) color center.  This defect is formed when a substitutional nitrogen is next to a lattice vacancy (missing carbon atom).  The center is usually found in one of two charge states, neutral and singly negatively charged, with the single negative charge state being the only one used for magnetometry.  The full details of the center are well reviewed in \citep{DOHERTY2013}.

The standard readout for NV is via some form of optically detected magnetic resonance (ODMR).  The NV electronic ground state is a spin triplet.  By illuminating off-resonantly (usually in the green) NV centers are optically polarised into the spin 0 state.  The degeneracy between spin 0 and spin $\pm 1$ states is lifted by a crystal field splitting around 2.88~GHz, and the Zeeman effect lifts the degeneracy between the $+1$ and $-1$ states, with increasing separation as a function of magnetic field.

The optical spin polarisation is due to the combination of spin conserving and spin non-conserving transitions, with the spin $\pm 1$ states more likely to de-excite via the spin-non-conserving singlet pathway.  This process is slower than direct spontaneous emission and occurs at 1040~nm, instead of the usual NV fluorescence spectrum from 630 to 750~nm.

A resonant radio-frequency field will drive transitions from spin 0 to one of the $\pm 1$ spin states.  The resonance will tend to mix populations, whereas the green excitation tends to polarise.  The combination of these effects means that when the RF field is in resonance, the fluorescence level drops, leading to a clear optical signal.

Lastly, we note that because of the crystallagraphic structure of diamond, there are four possible orientations for NV centers in the lattice.  Accordingly, the projection of an arbitrary magnetic field along each NV axis will be different, and so a single fluorescence diamond can be used as a vector magnetometer by reading out the magnetic field along each crystallographic axis.\citep{Schloss2018}.

Diamond magnetometry is a rapidly developing field, with many techniques designed to improve sensitivity including isotopic enrichment \citep{Balasubramanian2009}, portability through embedding in optical fibers \citep{Ruan2018,Bai2020,FILIPKOWSKI202210}, spin to charge readout \citep{Shields2015}, electrically detected magnetic resonance \citep{Bourgeois2015}, and laser threshold magnetometry \citep{Jeske2016,Dumeige:19,Hahl2022}.  As such, we have chosen to perform our subsequent analyses as a function of the magnetic field sensitivity, to better inform the sensitivities required for future platform development.

\section{INS aiding via map matching}
\label{sec2}
\begin{figure}[hpt]
  \centering
  \includegraphics[width=0.7\textwidth]{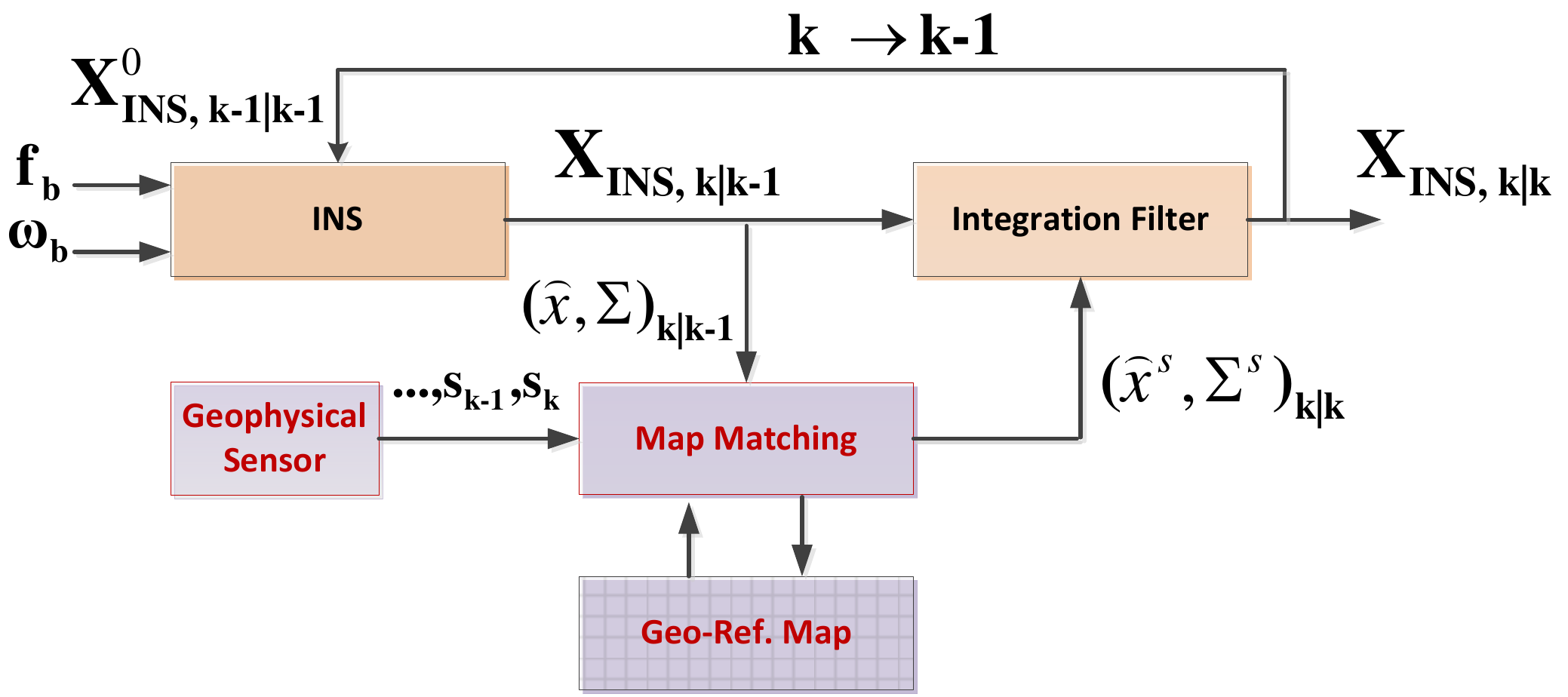}\\
  \caption{Illustration of a generic single recursion map matching aided Inertial Navigation System (INS). The map matching block takes the current kinematic state estimate $\hat{\bm{x}}$ and covariance $\bm{\Sigma}$ of the platform, and augments its position component by matching the geophysical measurement $\bm{s}_k$ to a geo-reference map. The prior kinematic state, $\hat{\bm{x}}$, is taken from the components of inertial navigation state $\bm{X}_{\mathrm{INS}}$, the latter is propagated using the accelerator and gyroscope measurements, respectively $\bm{f}_b$ and $\bm{\omega}_b$, based on the motion model on the surface of the earth in the INS block.} Finally, $\bm{X}_{\rm{INS}}$ denotes the full state vector from the INS, $\hat{\bm{x}}^s$ the estimated platform position and $\bm{\Sigma}^s$ the error covariance using $\bm{s}_k$ by the map matching.\label{fig-02}
\end{figure}

The process of an aided INS can be described by a recursive Bayesian filtering system, where the system prediction is given by the onboard INS, and system update is based on the measurements from external aiding sources. Fig.~\ref{fig-02} illustrates the block diagram of a generic aided INS with aiding from a map matching system. The INS is initialised from known parameters and at time $k$ propagates the navigation state $\bt{X}_{\mathrm{INS},k\mid k-1}$ based on the earth surface motion model and inertial measurements $(\bt{f}_b\,\bm{\omega}_b)$. The global position measurements, estimated from map matching, are assumed to be Gaussian distributed with mean $\hat{\bm{x}}^s$ - the estimated sensor location where $s$ is taken and covariance $\Sigma^s$, and are incorporated into the system via a widely used method to update the navigation state $\bt{X}_{\mathrm{INS},k\mid k}$.
For simplicity, we denote the navigation state at time $k$ as $\bt{X}_k \in \mathbb{R}^n$: this typically comprises the components of vehicle kinematic state (position and velocity) expressed in the geographical coordinates, vehicle attitude (roll, pitch and yaw), and inertial sensor bias terms. Based on inertial sensor measurements,  the propagation of the navigation state can be represented by
\begin{equation}\label{system}
\bt{X}_{k} = \bt{F}_{\rm INS}(\bt{X}_{k-1},\bm{\omega}_b,\bt{f}_b) + \bt{w}_k,
\end{equation}
where the function $\bt{F}_{\rm INS}(\cdot)$ signifies the mechanization of INS which involves the prior navigation state $\bt{X}_{k-1}$, the measurements of accelerometer $\bt{f}_b$ and gyroscope $\bm{\omega}_b$ at $k$, respectively \citep{Titterton2004a}, where $\bt{w} \sim \N(0,\bm{Q})$ accounts for system process noise including the errors from accelerometer and gyroscope\footnote{We have chosen this approach because it is widely used in the navigation literature \citep{Titterton2004a}, though the incorporation of accelerometer and gyroscope readings in this way is fraught with issues and leads to significant errors. This approach has the advantage of being generic; a  better method would require more realistic platform models.  We intend to discuss this issue in a future publication.}  $\N(\bt{a},\bm{B})$ signifies a Gaussian distribution with mean vector $\bt{a}$ and covariance matrix $\bm{B}$.

At each aiding update time $k$, the aiding position measurement is coupled into the navigation state via
\begin{equation}\label{meas1}
\bm{y}_k = \bt{H}\bt{X}_k + \bt{v}_k.
\end{equation}
where $\bt{H}$ is a constant matrix and $\bt{v} \sim \N(0,\bm{R})$ is a Gaussian zero-mean noise term modeling the measurement errors.

The INS aiding problem is to find the posterior density $p(\bt{X}_k\mid\bm{y}_{1:k})$ based on the sequence of measurements $\bm{y}_{1:k}$ from aiding sources.  Recursive Bayesian solution can be written as
\begin{equation}\label{s2-01}
p(\bt{X}_k\mid \bm{y}_{1:k}) = \frac{p(\bm{y}_k\mid \bt{X}_k)p(\bt{X}_k\mid \bm{y}_{1:k-1})}{\int p(\bm{y}_k\mid \bt{X}_k)p(\bt{X}_k\mid \bm{y}_{1:k-1}) d\bt{X}_k}.
\end{equation}
where the predicted density of the navigation state is given by
\begin{equation}\label{s2-011}
p(\bt{X}_k\mid \bm{y}_{1:k-1}) =  \int p(\bt{X}_k\mid \bt{X}_{k-1}) p(\bt{X}_{k-1}\mid\bm{y}_{1:k-1}) d\bt{X}_{k-1},
\end{equation}
and the navigation state prediction $\bt{X}_k\mid \bt{X}_{k-1}$  is a nonlinear process via the INS iterations. (\ref{s2-011}) is calculated by  (\ref{system}) as
\begin{equation}\label{s2-03}
p(\bt{X}_k\mid \bm{y}_{1:k-1}) \approx \N(\hat{\bt{X}}_{k\mid k-1},\, \bm{\Sigma}_{k\mid k-1})
\end{equation}
As shown in Fig.~\ref{fig-02}, at time $k$, the predicted density of the navigation state (\ref{s2-03}) is the output of the INS. For position aiding, the likelihood function $p(\bm{y}_k\mid \bt{X}_k)$ in (\ref{s2-01}) can be written as
\begin{equation}\label{s2-04}
p(\bm{y}_k\mid \bt{X}_k) = p(\bm{y}_k\mid \bm{x}_k).
\end{equation}
where $\bm{x}$ signifies the position components of $\bt{X}$. In the context of this work, $\bm{y}_k$ is the position measurement at $k$ to be obtained from map matching. The map matching problem is to find the posterior density of the vehicle position $\bm{x}_k$,  denoted by $p(\bm{x}_k\mid \bm{y}_{1:k})$,  using a TMI geo-referenced map $\mathcal{M}$ based on the magnetometer measurement sequence $\bm{y}_{1:k}$.

Similar to other geophysical data maps, the map matching localisation using TMI maps is facing the same challenges as we mentioned in Section~\ref{sec1}.
Next, we introduce a probabilistic data association method to address the map measurement ambiguity problem and carry out the ``curve fitting'' using a batch based probabilistic multiple hypothesis tracker.

\section{Probabilistic multiple hypothesis map matching}
\label{sec3}

The probabilistic multiple hypothesis map matching involves probabilistic data association (PDA) for data mapping from TMI signal domain to vehicle position domain, and a batch based multiple hypothesis tracking algorithm to iteratively optimise the estimated vehicle trajectory.

\subsection{Data mapping via probability data association}
Let $s_k$ represent the magnetometer measurement at time $k$. The measurement model is expressed as
\begin{equation}\label{meas}
s_k = s_k^o + \nu_k,
\end{equation}
where $s_k^o$ signifies the ground truth value of magnetic intensity measurement and $\nu_k$ is a noise term to cover the errors due to the imperfect of sensor. In this work, it is assumed to be a Gaussian distributed random variable, i.e., $ \nu \sim \N(0, \sigma^2) $.

According to (\ref{meas}), we need to consider a set of candidate measurements from a single measurement $s_k$, one of which is the true sensor measurement. Let $Z_m=\{\bm{z}_i, \, i = 1,\cdots, n\}$ denote the set of possible measurement locations on the map corresponding to $s_k$. In the meantime, we assume that at time $k$, the location of magnetometer, which takes the magnetic intensity measurement $s_k$, is a Gaussian random variable with mean $\bm{x}_k^s$ and covariance matrix $\Sigma^s_k$. Then, the location of true magnetic intensity measurement $\bm{z}_i, \, i = 1,\cdots, n$ should satisfy
\begin{equation}\label{s2-02}
 (\bm{z}_i-\bm{x}^s)(\bm{\Sigma}^s)^{-1}(\bm{z}_i-\bm{x}^s)' \leq \gamma,
\end{equation}
where $\gamma$ is a probability threshold, and its value determines an ellipsoid area on the data map containing the location of the magnetometer with a certain level of confidence. In this work, we refer to such an area as a search window.

\begin{figure}[htb!]
  \centering
  \includegraphics[width=0.7\textwidth]{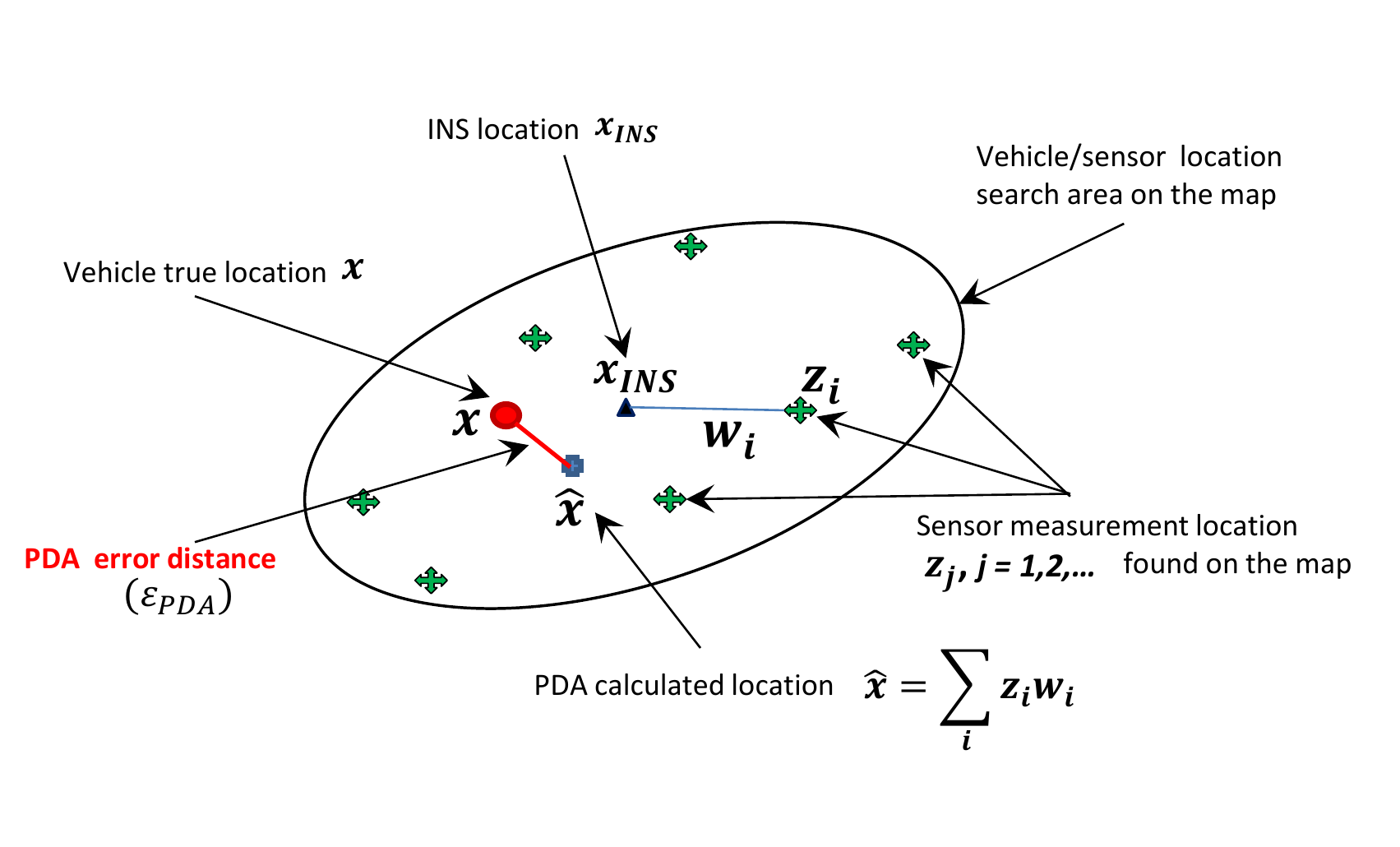}\\
  \caption{Illustration of the collection of candidate signal locations $\{\bm{z}^{s}_i,\,i = 1,2,\cdots,n\}$ obtained via (\ref{s2-01}) and (\ref{s2-02}) based on knowledge of predicted vehicle position $\bm{x}_{INS}$ from INS, and sensor noise level.}\label{fig-new1}
\end{figure}

PDA is a key method to determine where a magnetometer measurement originates on the geophysical data map in the presence of map measurement ambiguity.  Fig.~\ref{fig-new1}  illustrates the data PDA mapping process. As part of the integration filtering process shown in Fig.~\ref{fig-02}, the position components of the navigation state $\bm{X}_{INS}$, denoted by $\bm{x}_{INS}$, serve as the filter predicted position vector and will be updated using UKF if a position fix from an external source (e.g. GPS location, map matching position, etc.) is present. Note that if no external fix is present, the system shown in Fig.~\ref{fig-02} is a standard INS. $\bm{x}_{INS}$ and its covariance $\bm{\Sigma}$ describes a statistical bound (or uncertainty area) about where the signal $s$ is measured from the map. Therefore, a finite set of potential locations for signal $s$ on the map can be obtained via (\ref{s2-02}).

The probability weight of each candidate location $\bm{z}_i$ is proportional to the geometric distance between $\bm{z}_i$ and the window centre $\bm{x}^s$ (i.e., $\bm{x}_{INS}$).
The probability weight can be found as
\begin{equation}\label{s3-01}
 w_i = \frac{p(\bm{z}_i\mid\bm{x}^s)}{\sum_{j=1}^n p(\bm{z}_j\mid\bm{x}^s)},
\end{equation}
where $p(\bm{z}_i\mid\bm{x}^s) \sim \N\bigl(\bm{z}_i-\bm{x}^s,\,\bm{R}_i(\sigma)\bigr)$, and $\bm{R}_i(\sigma)$ is the associated variance which is a function of the signal noise variance, or in other words, signal-to-noise ratio (SNR).
Thus, the PDA solution for the map location on magnetic intensity measurement $s_k$ is the weighted mean:
\begin{equation}\label{s3-02}
\bar{\bm{z}} = \sum_{i=1}^n w_i \bm{z}_i.
\end{equation}
and the associated weighted variance:
\begin{equation}\label{s3-021}
\bar{\bm{R}} = \sum_{i = 1}^{n} w_i \left[\bm{R}_i(\sigma) + (\bm{z}_i-\bar{\bm{z}})(\bm{z}_i-\bar{\bm{z}})'\right].
\end{equation}

\subsection{Characterise map matching quality}

The performance of using PDA for magnetometer data mapping can be measured by the  PDA error distance $\varepsilon_{\rm PDA} $, which is defined as the Euclidian distance between the true magnetometer location\footnote{We assume that $\bm{x} = \bm{x}^s$ in the context of this work. Practically, we approximate the true magnetometer location by the position predicted by INS, i.e., $\bm{x}^s \approx \bm{x}_{INS}$.} and the location estimated via PDA.
\begin{equation}\label{s3-022}
\varepsilon_{\rm PDA} = \| \bt{\hat{x}}_{PDA}- \bm{x}^s\|.
\end{equation}
Note that (\ref{s3-022}) is a function of the magnetometer noise level $\sigma$, the location uncertainty of the magnetometer described by $\Sigma^s$ and resolution of the underlying geophysical data map.
\begin{figure}[bt!]
  \centering
  \includegraphics[width=0.7\textwidth]{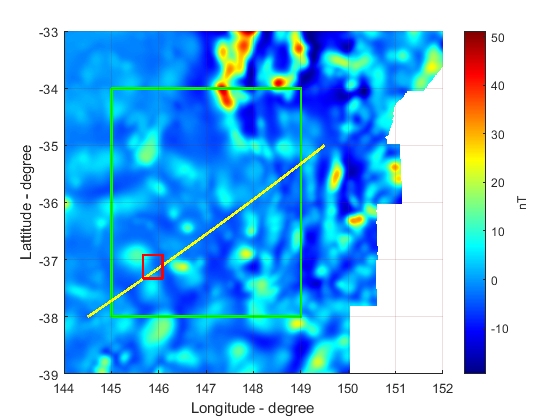}\\
  \caption{The total magnetic intensity map used in the simulation from \cite{TMI_online_map}. The platform travel trajectory (yellow line) is superimposed on the map. The green and red boxes illustrate the data areas considered in the simulations with results shown in Fig.~\ref{fig-61} and Fig.~\ref{fig-60}, respectively.}  \label{fig-51}
\end{figure}
Here, we use a simulation result to demonstrate this point. The total magnetic intensity map used in the simulation is downloaded from \cite{TMI_online_map}. As shown in Fig.~\ref{fig-51},  the actual data grid size is $85 \times 85$ metres. The simulation is carried out in the area surrounded by the green solid line rectangle. For every sensor noise level, 1000 samples are drawn randomly in the area, which are treated as the mean of sensor locations. The values of sensor location covariance $\Sigma^s$ and probability threshold $\gamma$ are chosen such that a search window approximately $6.8 \,\mbox{km}^2$ is formed for collecting candidate measurement locations. PDA error distances are then calculated as a function of sensor noise level $\sigma$ and map grid size.

\begin{figure}[bt!]
  \centering
\includegraphics[width=0.7\textwidth]{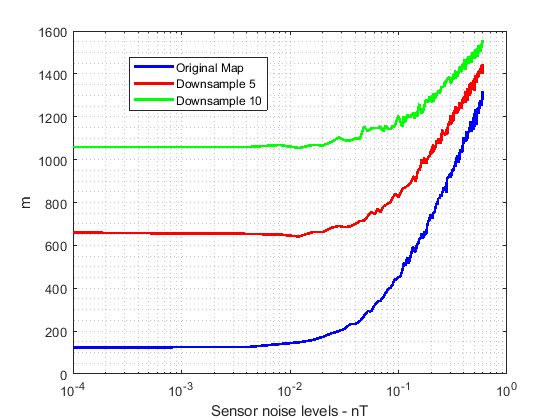}
\caption{Comparison of PDA error distances vs. sensor noise levels in the TMI map at original, 5 times and 10 times downsampled data grids.  The result is averaged from the results drawn in 1000 sample points in the green rectangle area shown in Fig.~\ref{fig-51}. For a given resolution, there is a plateau where increasing magnetometer sensitivity provides negligible improvement in position resolution.} \label{fig-61}
\end{figure}
Fig.~\ref{fig-61} shows the plot of PDA error distances versus sensor noise levels with original TMI map and the TMI maps of downsamped 5 and 10 times, respectively. The plots show that 1) map matching localisation error can be reduced if a magnetometer of low noise is used; 2) for this given map, about an accuracy of 150 metres may be achieved if a magnetometer of noise level 0.01 nT is used; 3) the map matching localisation error will increase if a TMI map of lower data grid is used, or if the magnetometer measurement noise level increases.

We use Map Feature Variability (MFV) as a measure of data variation sparsity of the geophysical data map. The MFV at a data point $i$ on a data map is defined as
\begin{equation}\label{s3-023}
\mathcal{C}_i = \frac{1}{n}\sum_{j}^n(s_{x_i} - s_{x_j})^2 \hspace{0.3 cm} \forall \, \bm{x}_j \in \mbox{search window},\, \bm{x}_j \neq \,\bm{x}_i.
\end{equation}

\begin{figure}[tb]
  \centering
  \includegraphics[width=0.9\columnwidth]{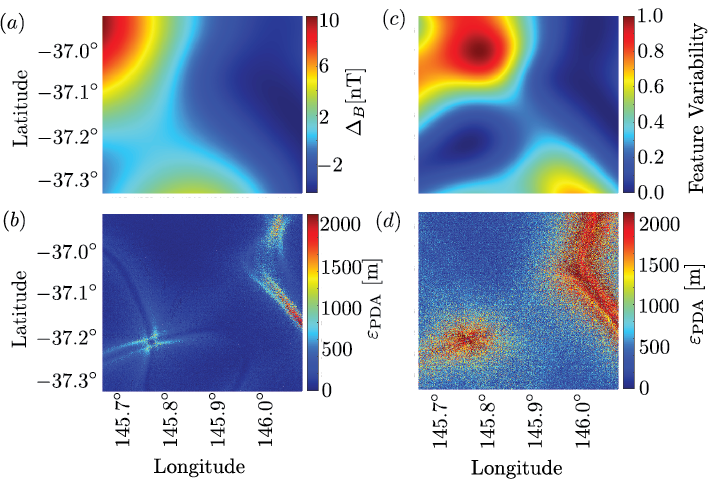}
  \caption{TMI map quality analysis in the red rectangle area shown in Fig.~\ref{fig-51}. (a) TMI map. (b) Map feature variability. (c) PDA error distance for $\sigma = 0.015$ nT. (d) PDA error distance for $\sigma = 0.15$ nT.  }\label{fig-60}
\end{figure}

Fig.~\ref{fig-60}(b) shows an example of the normalised map feature variability over the map area shown inside the red rectangle in Fig.~\ref{fig-51}. For reference, we also show the original TMI map in Fig.~\ref{fig-51}(a), and the PDA error distance maps for sensor noise levels $\sigma = 0.015$ nT and $0.15$ nT in Fig.~\ref{fig-51}(c) and (d), respectively.

From map matching localisation view point, MFV provides a local measure on map data variation sparsity. In a map matching based INS aiding, the value of $C_i^{-1}$ may be used to weight the estimated sensor location covariance to provide additional parameter that locally describes quality of the data map used.

\subsection{Map matching localisation}

We use a batch form Bayesian filter for map matching localisation from a sequence of magnetometer measurements. As shown in Fig.~\ref{fig-02}, the filter estimated posterior probability density of vehicle/sensor location $\N(\hat{\bm{x}}_k^s, \Sigma_k^s)$ based on a batch of $m$ magnetometer measurements  $\{s_{k-m+1},\cdots, s_{k}\}$ and a set of prior of the vehicle location $\{(\hat{\bm{x}}, \Sigma)_{k-m+1\mid k-m},\cdots, (\hat{\bm{x}}, \Sigma)_{k\mid k-1}\}$.

Two algorithms for such a filter may be used for this work, one is the probabilistic multiple hypothesis tracker based map matching (PMHT-MM) proposed in \citep{Wang2022} and the other is the Veterbi map matching algorithm proposed in \citep{Li2022}. Both algorithms work under a stochastic optimisation procedure which iteratively estimates the current vehicle location from a batch of measurements under vehicle dynamic constraints and are shown robust and accurate for INS navigation aiding with gravitational field data.

\section{Navigation experiment}
\label{sec4}

The simulation scenario is a constant velocity vehicle traveling along the surface of the earth at a fixed height of 100 m from $[-38^\circ, 144.5^\circ]$ to $[-35^\circ, 150^\circ]$ (i.e., from Melbourne area to Sydney area) and at a ground speed of 22 m/s.
The entire journey takes more than 3.6 hours and navigation is conducted by an onboard INS in GNSS denied environment.  The inertial sensors (both accelerometer and gyroscope) used in the INS are precision grade with errors specified in Table~\ref{table01}, and they take measurement at the frequency of 1 Hz and are assumed to be well calibrated before the gurney starts.
\begin{table}[ht]
\begin{center}
\caption{Bias and noise ranges of inertial sensors in the simulation according to \cite{2005:Christopher_Jekeli}.} \label{table01}
\begin{tabular}{@{}llll@{}}
  \toprule
  Sensor Grade & Sensor Type & Bias $b$ & White Noise $\sigma$ \\ \hline
 Precision (PC)  & Accel. horiz. & $2\times 10^{-6} m/s^2$ & $8\times 10^{-5} m/s^2/\sqrt{Hz}$ \\\cline{2-4}
   & Accel. Vert. & $2.5\times 10^{-8}m/s^2$ & $1.6\times 10^{-6} m/s^2/\sqrt{Hz}$ \\\cline{2-4}
   & Gyro. horiz.  & $2\times 10^{-5} deg/h$ & $1\times 10^{-3} deg/h/\sqrt{Hz}$ \\\cline{2-4}
   & Gyro. vert. & $1\times 10^{-3} deg/h$ & $3\times 10^{-2} deg/h/\sqrt{Hz}$ \\ \hline
\end{tabular}
\end{center}
\end{table}
\begin{figure}[tb!]
  \centering
  \includegraphics[width=0.78\textwidth]{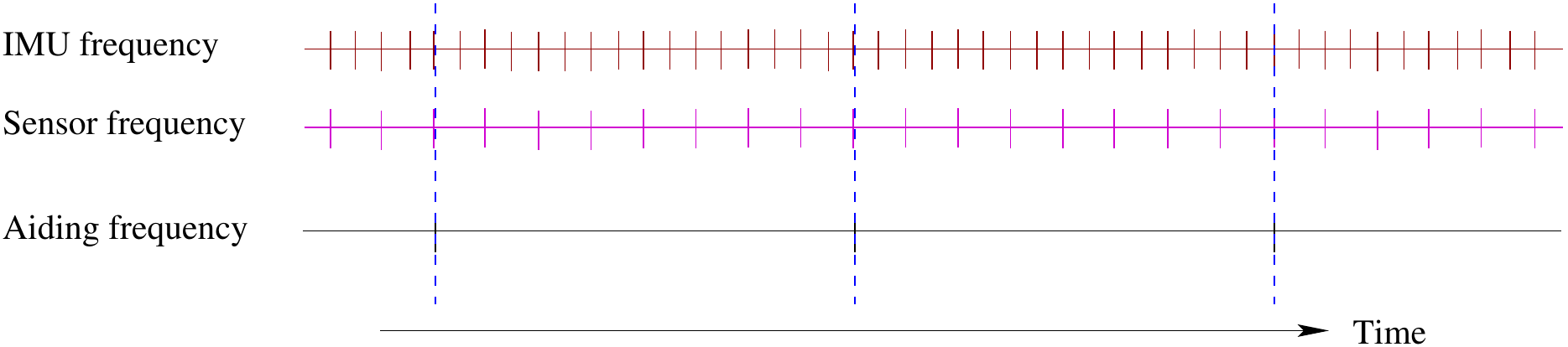}\\
  \caption{Illustration of the INS aiding timing in the simulation. We assume that measurement from inertial sensors is sampled at frequency 1 Hz and the magnetometer measurement is taken at an interval of every 10 seconds. Each time, the number of magnetometer measurements for PMHT-MM filtering are 30, which means that the actual aiding interval is 300 seconds.  }\label{fig-64}
\end{figure}
We assume that a low noise magnetometer is onboard to take magnetic intensity measurement at an interval of every 10 seconds. The PMHT-MM algorithm works with a batch of 30 magnetic intensity measurements at a time to estimate the vehicle current location and the estimated location is then used to update the navigation state via the UKF integration filter. This results in an aiding interval 300 seconds. Fig.~\ref{fig-64} illustrates the actual timing of data sampling in the simulation.

As we have shown in Fig.~\ref{fig-51}, a TMI map, which covers the area of southeast Australia, is obtained online from Geoscience Australia \citep{TMI_online_map} and is used in the simulation of map matching aided INS. The vehicle travel trajectory is shown by yellow solid line on the map.

\begin{figure}[bt!]
  \centering
  \includegraphics[width=0.7\textwidth]{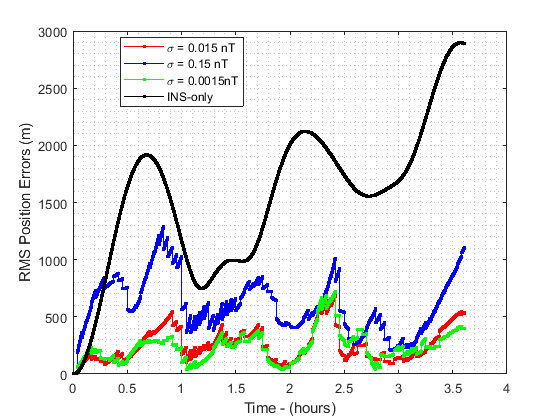}
  \caption{Comparison of RMS position errors of the aided INS for magnetometer measurement noise levels 0.015 nT (red) and 0.15 nT (blue), along with INS-only case (black). Note that the position error of INS-only case grows with time as a result of the inertial sensor bias, mainly from gyroscope, accumulating over time and its oscillation at a period approximately 1.4 hours is due to the Schuler pendulum effect occurring when navigation on a rotating frame \cite{Titterton2004a}.}   \label{fig-52}
\end{figure}

In this experiment, two noise levels for the magnetometer are considered: 1) $\sigma = 0.0015$ nT,  which is to model a magnetometer of very high precision; 2) $\sigma = 0.15$ nT, which represents the level of sensitivity of magnetometers that are commercially available. We plot the vehicle root-mean-squared (RMS) position errors in Fig.~\ref{fig-52} along with case of INS-only without aiding. The results were averaged from 100 Monte Carlo runs for each of cases.

The simulation results in Fig.~\ref{fig-52} show that
\begin{itemize}
\item when the noise level of magnetometer is 0,0015 nT, the INS with magnetometery aiding can achieve an average of RMS position error of 250 metres; the RMS position error will double if the noise level is 10 times larger at 0.15 nT.
\item we observed that the magnetometery INS aiding is robust with 100\% success rate with a batch length (i.e., the number of magnetometer measurements to be processed in a batch) 30 at each time in this simulation. If a lower batch length or a high noise level magnetometer is used,  the RMS position error will increases and the  magnetometery INS aiding will not be completely reliable.
\item the magnetometery INS aiding is able to remove the position drift, as indicated by the INS-only case, which is accumulated over time due to the imperfection of inertial sensors.
\end{itemize}

\section{Conclusions}
\label{sec5}

In this paper, a probabilistic method for map matching localisation based on magnetometery measurement and total magnetic intensity maps is described. We show that the method is able to effectively address the challenge issues associated with  map matching using geophysical maps and provides a mechanism of handling map measurement ambiguity and a way of evaluating the underlying quality. Furthermore, the effectiveness of the magnetometery map matching localisation is demonstrated using the simulation of removing position drift of an inertial navigation system, that arises in INS over a long duration, by the magnetometery aiding in the absence GNSS positioning. Simulation results using online maps verified the robustness and effectiveness of the proposed algorithm, particularly, the aiding precision will be getting better if a high sensitivity magnetometer is used.

\bmhead{Acknowledgments}
This work was funded by the Department of Defence Australia through the Next Generation Technology Fund.


\end{document}